\documentclass[11pt]{article}
\usepackage[utf8]{inputenc}
\usepackage{amsmath,amssymb,amsthm}
\usepackage{algorithm,algorithmic}
\usepackage[margin=1.1in]{geometry}
\usepackage{graphicx}
\input{macros.sty}
\usepackage{color}

\usepackage[textsize=small]{todonotes}
\newtheorem{remark}{Remark}
\providecommand{\keywords}[1]
{
  \small	
  \textbf{\textit{Keywords---}} #1
}

\title{{Randomized} Reduced Basis Methods for Parameterized Fractional Elliptic PDEs}
\author{Harbir Antil\thanks{Department of Mathematical Sciences and Center for Mathematics and Artificial Intelligence, George Mason University Email: hantil@gmu.edu} \and Arvind K.\ Saibaba\thanks{Department of Mathematics, North Carolina State University Email: asaibab@ncsu.edu}}

\begin{document}

\maketitle
\begin{abstract}
    This paper is interested in developing reduced order models (ROMs) for repeated simulation of fractional elliptic partial differential equations (PDEs) for multiple values of the parameters (e.g., diffusion coefficients or fractional exponent) governing these models. These problems arise in many applications including simulating Gaussian processes, {geophysical electromagnetics}. The approach uses the Kato integral formula to express the solution as an integral involving the solution of a parametrized elliptic PDE, which is discretized using finite elements in space and sinc quadrature for the fractional part. The offline stage of the ROM is accelerated using a solver for shifted linear systems, MPGMRES-Sh, and using a randomized approach for compressing the snapshot matrix. Our approach is both computational and memory efficient. Numerical experiments on a range of model problems, including an application to Gaussian processes, show the benefits of our approach. 
\end{abstract}


\section{Introduction}
We want to derive efficient methods for the parameterized fractional partial differential equation (PDE) 
\begin{equation}\label{eqn:paramfpde}
    (\mc{A}(\bfx;\bfmu))^{\alpha}y(\bfx;\bfmu) = b(\bfx;\bfmu), \qquad \bfx \in \Omega
\end{equation} 
where $\mc{A}$ is an elliptic differential operator, parameterized by the set of parameters $\bfmu \in \mc{P} \subset \R^{p}$. Here $0 < \alpha < 1$ is the fractional exponent. Throughout this paper, unless otherwise specified, we consider the boundary conditions to be homogeneous Neumann, but the approach can be modified slightly to accommodate other boundary conditions as well. {We assume that the fractional power is defined in the spectral sense.} Solving the fractional elliptic PDE is a computationally challenging problem in itself. {Indeed, \eqref{eqn:paramfpde} is nonlocal and traditional `local basis' approaches such as finite element method require dealing with dense numerical linear algebra. Therefore, the efficiency of traditional iterative methods becomes questionable.
}

In this paper, we want to address the more challenging problem of efficiently solving~\eqref{eqn:paramfpde} for multiple instances of the parameter $\bfmu$ and the exponent $\alpha$. This is accomplished by the use of the reduced basis approach. These kinds of problems arise in {imaging \cite{HAntil_SBartels_2017a}, geophysical electromagnetics \cite{CJWeiss_BGVBWaanders_HAntil_2020a}, fractional Laplacian regularizers in optimization \cite{HAntil_ZDi_RKhatri_2020a}, harmonic maps \cite{HAntil_SBartels_ASchikorra_2022a}}. In~\cite{trillos2017bayesian}, the authors investigate the well-posedness of a Bayesian inverse problem with a fractional elliptic PDE forward model. To estimate the diffusion coefficient from data, multiple instances of the fractional PDE have to be {solved. A similar situation occurs in optimal control problems where the fractional exponent is considered as a random variable \cite{HAntil_DPKouri_JPfefferer_2021a}, or optimization problems to identify the fractional exponent \cite{HAntil_EOtarola_AJSalgado_2018b,sprekels2016new}. Another application of interest is the stochastic PDE (SPDE) approach to Gaussian processes, which is explored further in Section~\ref{ssec:gp}. In view of all these applications, there is a clear need for an efficient reduced order model (ROM).}

In many applications of interest, high-fidelity models are governed by a single PDE or systems of PDEs that are expensive to simulate. Reduced order models refer to a collection of techniques for efficiently reducing the complexity of high-fidelity models (also known as full-order models or FOMs) with perhaps only a small sacrifice in accuracy compared to the FOMs.  ROMs are especially beneficial for problems in which the full-order model has to be simulated for several parameter conditions as the ``inner loop'' for a variety of applications such as optimization, control, inversion, uncertainty quantification, design, or digital twins. This is known as parametric ROM and this is the setting in the present paper. 

While there are several possible techniques for ROMs, we follow the so-called projection-based model reduction approach, in which there is a clear separation into two stages. The first stage is called the {\em offline stage} in which the FOM is simulated under a wide range of ``representative'' conditions of the parameters to collect snapshots of the system. This is then followed by a compression step to extract a basis that is used to then project the FOM onto a lower-dimensional subspace yielding a ROM. In the {\em online stage}, the lower-dimensional ROM is then simulated for a given parameter. The computational tradeoff comes from the fact that the offline costs may be large but upfront, and acceptable in certain situations, so long as the online costs are small and the accuracy is good. 

While ROM for elliptic PDE-based problems is a well-established field, ROM for fractional elliptic PDEs is relatively less well-studied and is the main focus of this paper.  We briefly review the literature in this area and summarize the main contributions of this paper.  

\paragraph{Literature review}
{A recent paper~\cite{antil2019reduced} considers a reduced basis approach for the spectral definition when realized using the Stinga-Torrea (Caffarelli-Silvestre) extension \cite{LCaffarelli_LSilvestre_2007a,PRStinga_JLTorrea_2010a}. The articles~\cite{danczul2021reduced,danczul2022reduced} employ a reduced basis approach by interpolating operator norms. For another reduced ordering modeling approach for problems with nonlocal integral kernels, we refer to \cite{MR3612777}. The paper~\cite{fu2018pod} considers ROMs for time-fractional PDEs. The recent review paper~\cite{d2020numerical} reviews numerical methods, including reduced order models, for a different definition of the nonlocal problem.
}

The following papers~\cite{dinh2019model,bonito2020reduced} are closest to our setting, in which  $\mc{A}(\bfx;\bfmu)y(\bfx;\bfmu) = (-\Delta y)(\bfx)$ is the standard Laplacian and $g(\bfx;\bfmu) = g(\bfx)$. The reduced order model is constructed over the exponent $\alpha$, but the parametric dependence of the PDE on other parameters $\bfmu$ is not considered. {Notice that, these papers build on the traditional deterministic reduced order method. Furthermore, in our setting, we use randomized approaches, which offer additional benefits as highlighted below. }
In recent work~\cite{antil2022efficient}, we have used the MPGMRES-Sh approach~\cite{bakhos2017multipreconditioned} for efficiently applying the fractional operator for a single parameter value $\bfmu$. It was observed in this paper that the same basis can be used for multiple values of the exponent, which is important in our approach as well.

We also briefly mention that the idea of using randomization in ROMs has  previously been explored; see, e.g., ~\cite{saibaba2020randomized} and references within. The randomized compression technique we use in this paper is applicable to the streaming setting, in which several columns of the snapshot matrix are added sequentially.  

\paragraph{Contributions} In this paper, we derive reduced order models for solving~\eqref{eqn:paramfpde}. Our approach is based on the following features:
\begin{enumerate}
    \item We use the spectral definition of the fractional elliptic operator and use the Kato integral formula to express the solution of~\eqref{eqn:parampde} as an integral involving the solution of a parametric elliptic PDE with an additional parameter to compute the integral. The operator is then discretized using standard Lagrange finite elements, and the integral is discretized using the sinc quadrature approach. 
    \item We use a reduced basis approach for efficiently computing the solutions of the parametric elliptic PDE. It contains an offline stage in which a reduced basis is generated, and in the online stage, a projected system involving the reduced basis is solved. 
    \item The construction of the offline stage exploits the structure of the problem. We use an efficient iterative method MPGMRES-Sh to accelerate the convergence of shifted systems; this approach drastically reduces the number of PDE solves needed for generating the basis.  Furthermore, we use randomized solvers to mitigate the cost of storing and compressing the snapshots, and an efficient method involving generalized eigendecomposition to accelerate the costs of the online stage.
\end{enumerate}
{As a part of future work, we plan to explore rigorous a posteriori error analysis \cite{JSHesthaven_GRozza_BStamm_2016a,AQuarteroni_AManzoni_FNegri_2016a} in the context of our randomized ROM approach.}

\section{Background}
In this section, we discuss the integral approach to solving the fractional elliptic PDE. We also discuss the ingredients required for the ROM, namely the iterative solver MPGMRES-Sh (Section~\ref{ssec:mpgmressh}) and the randomized approach we use for compressing the snapshots (Section~\ref{ssec:rand}).
\subsection{Integral approach}\label{ssec:setup}
 Here we follow the approach in~\cite[Section 2.2]{dinh2019model} that combines an integral formulation for the inverse fractional operator with a parameterized elliptic PDE. Using the Kato integral formula~\cite[Theorem 2]{kato1960note}, with a simplification $\lambda=0$, the solution to~\eqref{eqn:paramfpde} can be expressed as  
\[ y(\bfx;\bfmu) = \frac{\sin(\alpha \pi)}{\pi}\int_{-\infty}^\infty e^{(1-\alpha)z} (e^z + \mc{A}(\bfx;\bfmu))^{-1}b(\bfx;\bfmu)dz. \]
After setting  $\bfrho = (\bfmu, z)$, we get the parametric PDE
\begin{equation}\label{eqn:parampde}
    (e^z +  \mc{A}(\bfx;\bfmu))u(\bfx;\bfrho) = b(\bfx;\bfmu) \qquad \bfx \in \Omega,
\end{equation}
with appropriate boundary conditions, and arrive at the following integral formula {for $y$} 
\begin{equation}\label{eqn:integralw} y(\bfx;\bfmu) = \frac{\sin(\alpha \pi)}{\pi}\int_{-\infty}^\infty e^{(1-\alpha)z}u(\bfx;\bfrho) dz.\end{equation}
This integral approach is the foundation for the ROM that we construct. Following~\cite{bonito2015numerical}, we use finite element discretization in space and sinc quadrature for the integral, which we now discuss.
\paragraph{Discretization} We pick a finite-dimensional basis for the subspace $\mc{V}_h \subset L^2(\Omega)$ denoted by $\{\phi_j\}_{j=1}^{N_h}$ so that $\textrm{dim}(\mc{V}_h) = N_h$ and any function $v_h \in V_h$ can be represented by its expansion in the basis $v_h= \sum_{j=1}^{N_h}v_{h,j}\phi_j$. We denote the vector of its coefficients as $\bfv = \bmat{v_{h,1} & \dots & v_{h,N_h}}^\top$. For the sinc quadrature to discretize the integral representation~\eqref{eqn:integralw}, let $h$ represent the mesh parameter and define \begin{equation}\label{eqn:zlims}Z_+ = \lceil \frac{\pi^2}{4\alpha \zeta^2 }\rceil,\qquad Z_- = \lceil \frac{\pi^2}{4(1-\alpha) \zeta^2 }\rceil, \end{equation}
where $\zeta = 1/\log(1/h)$. We also define the quadrature weights and nodes as 
\begin{equation}\label{eqn:quadwts} w_j =  \frac{\zeta\sin(\alpha\pi)}{\pi} e^{(1-\alpha)z_j} \qquad z_j = \zeta j \qquad -Z_- \leq j \leq Z_+.  \end{equation}

Let $a(u,v; \bfmu)$ denote the weak form of the operator $\mc{A}(\bfx;\bfmu)$ with the corresponding boundary conditions; similarly, let $g(v;\bfmu)$ be the weak form of $b(\bfx;\bfmu)$. We seek a solution $u_{h,j} \in V_h$ such that 
\begin{equation}\label{eqn:discweak}
    a(u_{h,j}(\bfx;\bfmu) ,v_h(\bfx); \bfmu) + e^{z_j}\int_\Omega  u_{h,j}(\bfx; \bfmu) v_h(\bfx) d\bfx = g(v_h;\bfmu), \qquad \text{for all } v_h \in V_h 
\end{equation}
and for $-Z_- \leq j \leq Z_+$. Having obtained the solutions $u_{h,j}$, the approximate solution $y_h$ to~\eqref{eqn:integralw} can be obtained using the sinc-quadrature formula as $y_h(\bfx;\bfmu) =\sum_{j=-Z_-}^{Z_+} w_ju_{h,j}(\bfx;\bfmu)$. However, as can be seen, computing this quadrature formula requires the solution of many PDEs of the form~\eqref{eqn:discweak}. This can be computationally intensive if this process has to be repeated for multiple values of the exponent $\alpha$ and the parameter $\bfmu$.  This motivates the need for a reduced-order approach. We first discuss an efficient method for solving the linear systems arising from the system~\eqref{eqn:discweak}.

{The reference \cite{bonito2015numerical} establishes exponential convergence
in quadrature and also derives appropriate error estimates with respect to spatial
discretization. Provided that the so-called Kolmogorov $n$-width in the parameter
space decays exponentially, which is indeed the case when the solution depends 
smoothly on parameters (cf.~\cite{APinkus_1985}), one can establish that reduced 
basis error also decays exponentially. In the case that the reduced bases are 
generated using a greedy method, we refer to \cite{PBinev_ACohen_RDevore_2012a} 
for such results.}

\subsection{MPGMRES-Sh} \label{ssec:mpgmressh}
MPGMRES-Sh~\cite{bakhos2017multipreconditioned} is an iterative method for solving a sequence of shifted systems of the form 
\begin{equation}\label{eqn:shift2} (\bfC_1 + \sigma_j\bfC_2)\bfx_j = \bfb \qquad 1 \leq j \leq N_\sigma,\end{equation}
for a set of shifts $\{\sigma_j\}_{j=1}^{N_\sigma}$ and matrices $\bfC_1,\bfC_2 \in \R^{n\times n}$. 
The main idea of this approach is to use shift-and-invert preconditioners of the form $\{\bfC_1+\tau_j\bfC_2\}_{j=1}^{n_p}$, where $n_p$ is the number of preconditioners used. A simplified version of the algorithm is given in Algorithm~\ref{alg:mpgmressh}. 

\begin{algorithm}[!ht]
    \begin{algorithmic}[1]
        \REQUIRE Matrices $\bfC_1,\bfC_2$, Shifts $\{\sigma_j\}$ for $1 \leq j \leq {N_\sigma}$, Vector $\bfb$, Preconditioners $\bfP_j = \bfC_1 + \tau_j\bfC_2$, $1 \leq j \leq n_p$
        \STATE Compute $\beta = \|\bfb\|_2$ and $\widehat\bfV^{(1)} = \bfb/\beta$
        \FOR {$k=1,\dots ,$ until convergence}
        \STATE $\widehat{\bfv}_k = \widehat\bfV^{(k)}\bfe_{n_p}$ and $\bfW^{(k)} = \bmat{\bfP_1^{-1}\widehat{\bfv}_k & \dots & \bfP_{n_p}^{-1}\widehat{\bfv}_k}  $
        \STATE $\widehat{\bfW} = \bfC_2 \bfW^{(k)}$
        \FOR {$j=1,\dots,k$}
            \STATE $\bfH^{(j,k)} = (\widehat\bfV^{(j)})^\top\widehat\bfW$
            \STATE $\widehat{\bfW} = \widehat{\bfW} - \bfV^{(j)}\bfH^{(j,k)}$
        \ENDFOR
        \STATE $\widehat\bfW = \bfQ\bfR$ \COMMENT{Thin QR factorization}
        \STATE Construct $\widehat\bfV^{(k+1)} = \bfQ$ and $\bfH^{(k+1,k)} = \bfR$
        \STATE Check for convergence for all $\{\sigma_j\}_{j=1}^{N_\sigma}$ and terminate, if converged, at iteration $k_m$. 
        \ENDFOR 
         \STATE Consolidate matrices $\bfZ_{k_m} = \bmat{\bfW^{(1)} & \dots & \bfW^{(k_m)}}$
        \RETURN Basis $\bfZ_{k_m}$, where $k_m$ is the number of iterations upon termination
    \end{algorithmic}
    \caption{Simplified version of MPGMRES-Sh}
    \label{alg:mpgmressh}
\end{algorithm}
The idea behind this approach is to form a basis for the search space $\bfZ_{k_m} \in \R^{n\times k_mn_p} $ where $k_m$ is the number of MPGMRES-Sh iterations and at each iteration, we increment the dimension of the search space by $n_p$. Sometimes the basis fails to be incremented by $n_p$ but this can be handled using rank-revealing QR factorizations; we do not consider this in the present work but the reader is referred to~\cite{bakhos2017multipreconditioned}. This search space is common to all the shifted systems and is used to solve a projected subproblem for each shift $\sigma$ of size $(k_m n_p) \times (k_m n_p)$ to obtain a solution $\bfp(\sigma) \in \R^{(k_m n_p) \times 1}$. In a bit more detail, the approximate solution is computed as $\bfx_{k_m}(\sigma) = \bfZ_{k_m}\bfp_{k_m}(\sigma)$, where $\bfp_{k_m}$ solves the optimization problem
\[ \min_{\bfp \in  \R^{(k_mn_p)\times 1}} \|\bfb - (\bfC_1 + \sigma_j\bfC_2)\bfZ_{k_m}\bfp \|_2^2,  \]
for each shift $\sigma_j$, $1\leq j \leq N_\sigma$. In the application of interest in this present paper, we will focus on computing the basis $\bfZ_{k_m}$ rather than the solutions $\bfx(\sigma)$ as we now explain in this remark.

\begin{remark}\label{rem:mpgmpres}
Since $\bfZ_{k_m}$ is the search space corresponding to the shifted systems~\eqref{eqn:shift2}, we can write $\bfx_j \approx \bfZ_{k_m}\bfp_{k_m}(\sigma_j)$. Therefore, the matrix
\[ \bfX = \bmat{\bfx_1 & \dots & \bfx_{N_\sigma}} \approx \bfZ_{k_m} \bmat{\bfp_{k_m}(\sigma_1) & \dots & \bfp_{k_m}(\sigma_{N_\sigma})}. \]
That is, the range of $\bfX$ is approximately captured by the matrix $\bfZ_{k_m}$. Since the leading left singular vectors approximate the range of $\bfX$, it is sufficient to store $\bfZ_{k_m}$.  This key insight is important to reduce computational costs in the offline stage.
\end{remark}

 The number of iterations is determined when all the relative residual for the shifted systems is smaller than $10^{-8}$.  For more implementation details and theoretical justifications, the reader is referred to~\cite{bakhos2017multipreconditioned}. If the number of MPGMRES-Sh iterations, $k_m$, and the number of preconditioners, $n_p$, are both small, then the cost of solving the projected problem across multiple shifts is small and the cost of MPGMRES-Sh is dominated by the cost of forming the basis for the search space. In our applications, the number of preconditioners is $3$ and the number of iterations is smaller than $30$ yielding a basis $\bfZ_{k_m}$ of size at most $90$. We now discuss the second main ingredient in our approach, the randomized low-rank approximations.

\subsection{Randomized low-rank approximation}\label{ssec:rand}
A reduced-order model requires compression of the matrix containing snapshots of the state. In this paper, we pursue a randomized approach since it has reduced computational and storage requirements compared to the full SVD. Furthermore, as we will explain, this approach is beneficial in our setting  where the matrix is available in a streaming fashion.

Suppose we wanted to compute a low-rank approximation to the matrix $\bfS \in \R^{m\times n}$. We follow the approach from~\cite{tropp2017practical}.  We draw two random matrices $\bfOmega \in \R^{n\times \ell_1}$ and $\bfPsi \in \R^{m\times \ell_2}$; we take them to be standard Gaussian random matrices, which means that the entries are independent and drawn from the Gaussian distribution with zero mean and unit variance. We form the sketches $\bfY_1 = \bfS \bfOmega$ and $\bfY_2 = \bfS^\top \bfPsi$, which approximate the column and row spaces of $\bfS$ respectively. Then compute the thin-QR factorization $\bfY_1 = \bfQ\bfR$ and discard $\bfR$. The low-rank approximation of $\bfS$ can be computed as 
\[ \bfS \approx \bfS_{\rm rand}\equiv \bfQ (\bfPsi^\top \bfQ)^\dagger \bfY_2^\top,\]
where $^\dagger$ is the Moore-Penrose inverse. The low-rank approximation can be converted into the SVD format to extract singular vectors with relatively less additional work. If a rank-$r$ approximation of the matrix $\bfS$ is desired, assuming that the matrix is stored in dense form, the computational cost of this algorithm is  
$\mc{O}(rmn + r^2 (m+n) + r^3)$ flops.  The choice of the sketch sizes $\ell_1$ and $\ell_2$ is discussed in~\cite[Section 4.4]{tropp2017practical} which we follow here as well. We choose $\ell_1 = 2r +1 $ and $\ell_2 = 2\ell_1 + 1$. Based on this sketch size, the error in expectation of the low-rank approximation satisfies~\cite[Theorem 4.3]{tropp2017practical} 
\[ \mathbb{E}\|\bfS-\bfS_{\rm rand}\|_F^2 \leq 4 \left( \sum_{j >r }\sigma_j^2\right), \]
where $\{\sigma_j\}_{j=1}^{\min\{m,n\}}$ are the singular values of the snapshot matrix $\bfS$. The main message is that the error in the low-rank approximation due to the randomization is, in expectation, within a small factor of the ``best'' low-rank approximation (note that this approximation need not be unique). In particular, if the singular values of the snapshot matrix decay rapidly then the error in the low-rank approximation is very small. 

Besides lower computational costs compared to the full SVD algorithm (i.e., $\mc{O}(mn\min\{m,n\})$), the proposed approach is amenable to the streaming setting where the matrix is updated sequentially. That is, if the matrix $\bfS$ is of the form 
\[ \bfS = \bfS_1 + \dots + \bfS_M,\]
where the summands $\bfS_j$ are available sequentially. As an example that is relevant to our application, consider the case, where the columns (or a block of columns) of a matrix are computed sequentially.   It is easy to see that the sketches $\bfY_1$ and $\bfY_2$ can also be updated sequentially. Furthermore, once the sketches have been computed, the appropriate columns can be discarded since all the randomized algorithm only requires the sketches to compute the low-rank approximation. These insights  will be relevant to the proposed algorithm resulting in an algorithm with lower storage and memory costs.

\section{Reduced basis approaches}
In this section, we introduce the weak form and discretization of the parametric PDE that is the foundation of our ROM (Section~\ref{eqn:parampde}). Then we present two approaches for ROM: a na\"ive approach (Section~\ref{ssec:naive}), and an efficient approach (Section~\ref{ssec:prop}). We conclude this section with a discussion on computational costs. 

\subsection{Weak form and discretization}\label{ssec:parampde} We aim to solve the parameterized PDE~\eqref{eqn:parampde} for multiple instances of the parameters $\bfrho$. To this end, with an eye towards numerical approximation using finite elements, we consider the weak form of the problem, which reads as: given $\bfrho \in \mc{P} \times\R_+ $ find {$u(\bfrho) \in V$} such that
\[ a(u(\bfrho), v; \bfrho)  = g(v;\bfrho) \qquad \forall v \in V,  \]
where the exact form of $a$ depends on the form of $\mc{A}$ and the boundary conditions. We assume that the weak form has an affine parameter dependence as 
\begin{equation} \label{eqn:affine_a}
    a(u, v; \bfrho) := \sum_{t=1}^{n_a} f_t^a(\bfmu) a_t(u, v) +  e^z\int_\Omega uv d\bfx, 
\end{equation}
where $a_t: V \times V \rightarrow \R$ are independent of the parameters
and the right-hand side also satisfies an affine parametric dependence
\begin{equation} \label{eqn:affine_g}
    g( v; \bfrho) := \sum_{t=1}^{n_g} f_t^g(\bfmu) g_t(v).
\end{equation}
If the problem does not admit an affine parametric dependence, one can use techniques such as the empirical interpolation method~\cite[Chapter 10]{AQuarteroni_AManzoni_FNegri_2016a}. 

 We proceed to discuss the matrix representation of this parametric PDE that will be starting point of the ROMs. As before, let $V_h \subset V$ be a finite-dimensional subspace of $V$ with $\dim(V_h) = N_h$ with a basis $\{\phi_1,\dots,\phi_{N_h}\}$.  We denote the matrix representation of the bilinear forms as 
\begin{equation}\label{eqn:bilinear_mat}
    [\bfA_t]_{ij} = a_t(\phi_i,\phi_j) \qquad 1 \leq i,j \leq N_h, \quad 1 \leq t \leq n_a.
\end{equation}
We denote by $\bfM$, the mass matrix $\bfM \in \R^{N_h\times N_h}$  with entries $[\bfM]_{ij} = \int_\Omega \phi_i\phi_j d\bfx$ for $1 \leq i,j \leq N_h$. We also define the vectors $ \bfg_t$ with entries $[ \bfg_t]_i = g_t(\phi_i)$ for $1 \leq i \leq N_h$, and $1 \leq t \leq n_g$. In summary, the discretized parametric system we aim to solve is 
\begin{equation}
    \left(\sum_{t=1}^{n_a} f_t^a({\bfmu}) {\bfA}_t + e^z \bfM \right)\bfu = \sum_{t=1}^{n_g} f_t^g(\bfmu) \bfg_t, 
\end{equation}
for the set of parameters $\bfrho = (\bfmu,z)$. We now discuss ROMs for the parametric PDEs.  
\subsection{A na\"ive first approach}\label{ssec:naive} Before we present our efficient approach, we first discuss a na\"ive approach that will be the basis for comparison with our proposed approach.

We first draw a set of random instances of the parameters $\mc{M} = \{ \bfmu^{(j)}\} $ for $1\leq j \leq M$; strategies for sampling are discussed in~\cite[Section 6.6-6.7]{AQuarteroni_AManzoni_FNegri_2016a}. Similarly, using the definitions in Section~\ref{ssec:setup} we set the  parameters $\mc{Z} = \{z_j\}$ for $-Z_- \leq j \leq Z_+$ with $Z = Z_-+Z_++1$ corresponding to the setting $\alpha = 0.5$. This is motivated by the observation in~\cite[Section 5.1]{antil2022efficient} that the exponents $e^{z_k}$ are completely independent of the exponent, except for the limits $Z_-$ and $Z_+$.   As we will see in the numerical experiments, it is sufficient to generate the projection basis corresponding to a single exponent $\alpha$, which is then applicable to a wider range of the exponent.

Thus, we now have the set of parameters $\mc{M} \times \mc{Z}$ for $\bfrho$. 
For each parameter value $\bfmu^{(j)}$ we first construct the matrices and vectors
\begin{equation}\label{eqn:matvecs} \bfK_j =\sum_{t=1}^{n_a} f_t^a({\bfmu}_j) \bfA_t, \qquad \bff_j = \sum_{t=1}^{n_g} f_t^g({\bfmu}_j) \bfg_t, \quad 1 \leq j \leq M.\end{equation}
Then we solve the sequence of shifted linear equations 
\begin{equation}\label{eqn:shiftedsys} (\bfK_j + e^{z_k}\bfM)\bfu_{j,k} = \bff_j \qquad 1 \leq j \leq M, \quad -Z_- \leq k \leq Z_+, \end{equation}
to generate the snapshot matrix 
\begin{equation}\label{eqn:snapshot}
    \bfS = \bmat{\bfU^{(1)} & \dots  & \bfU^{(M)}  } \in \R^{N_h \times (MZ)}, \qquad \bfU^{(j)} =\bmat{ \bfu_{1,-Z_-}  & \dots & \bfu_{j,Z_+}} \qquad 1 \leq j \leq M.
\end{equation} 
The next step is to compute the dominant left singular vectors of $\bfS$ denoted $\bfV \in \R^{N_h \times K}$. We form the matrices $\widehat\bfA_t = \bfV^\top\bfA_t\bfV $ for $1\leq t \leq n_a$, $\widehat\bfM = \bfV^\top\bfM\bfV$, and vectors $\widehat\bfg_t = \bfV^\top\bfg_t$ for $1 \leq t \leq n_g$. 

This basis $\bfV$ can be used to approximate the solution space in the online stage. Given a new parameter value $\overline{\bfmu}$, and an exponent $\alpha \in (0,1)$, we solve the projected system (see~\eqref{eqn:zlims} for definitions of $Z_+,Z_-$)
\begin{equation}\label{eqn:projsys}
    \left(\sum_{t=1}^{n_a}f_t^a(\overline{\bfmu}) \widehat{\bfA}_t + e^{z_k} \widehat{\bfM} \right)\widehat\bfu_k = \left( \sum_{t=1}^{n_g} f_t^g(\overline{\bfmu})\widehat{\bfg}_t \right) \qquad -Z_- \leq k \leq Z_+.
\end{equation} 
The approximation solution is $\widehat{\bfy} =  \bfV \left( \sum_{k= -Z_-}^{Z_+}w_k\widehat\bfu_k\right)$. The details of this approach are given in Algorithm~\ref{alg:romonline}.

\begin{algorithm}[!ht]
\begin{algorithmic}[1]
\REQUIRE Sample parameter $\overline{\bfmu}$, Precomputed matrices $\{\widehat\bfA_t\}_{t=1}^{n_a}$ and vectors $\{\bfg_t\}_{t=1}^{n_g}$, exponent $\alpha$
\STATE Compute quadrature nodes $\{z_k\}_{-Z_-}^{Z_+}$ and weights $\{w_k\}_{-Z_-}^{Z_+}$ corresponding to the exponent $\alpha$; see~\eqref{eqn:zlims} and~\eqref{eqn:quadwts}.
\STATE Form and solve the projected system~\eqref{eqn:projsys}
\STATE Evaluate  $\widehat\bfw = \bfV \left(\sum_{k=-Z_-}^{Z_+} w_k \widehat\bfu_k \right)$.
\RETURN Approximate solution $\widehat\bfw$
\end{algorithmic}
\caption{Reduced basis approach for fractional elliptic PDE: Online stage}
\label{alg:romonline}
\end{algorithm}

However, the cost of generating and compressing the snapshot matrix can be prohibitively expensive in practice as we now explain. This is because  $MZ$ linear systems need to be solved in total to compute the snapshot matrix. Note that $M$ can range from $10-100$ whereas $Z$ can range from $100-1000$, depending on the value of the exponent $\alpha$ and the discretization of the system. Furthermore, the size of  $\bfS$ can be large so the cost of compressing this system is large; the cost of SVD is $\mc{O}(N_hMZ \min\{N_h,MZ\})$. Even the cost of storing the matrix may be expensive in practice. To address these issues we propose a new approach for efficiently generating and compressing the basis.  

\subsection{Proposed approach}\label{ssec:prop} There are two components to the proposed approach in the offline stage: the use of MPGMRES-Sh to solve the shifted systems and randomization to compress the snapshot matrix. In the online stage, we use a generalized eigendecomposition to speed up computational costs. 

\paragraph{Offline stage} First, we solve the systems~\eqref{eqn:shiftedsys} for each value of $\bfmu^{(j)}$ using the MPGMRES-Sh approach (Section~\ref{ssec:mpgmressh}). The algorithm is implemented with $\bfC_1 = \bfM$ and $\bfC_2 = \bfA_j$ and shifts $\sigma_k = 1/\exp(z_k)$ for $-Z_-\leq k \leq Z_+$. We get a basis $\bfZ^{(j)} \in \R^{N_h \times (n_m)}$ where $n_m$ is the size of the MPGMRES-Sh basis (obtained as $n_m = k_mn_p$, i.e., the number of MPGMRES-Sh iterations $k_m$ times the number of preconditioners $n_p$). For simplicity, we have assumed that the MPGMRES-Sh basis size is the same for every $j$ but in practice, it may be different. We have also dropped the dependence of the basis on the number of iterations $k_m$.  
In all the experiments, for MPGMRES-Sh we take three preconditioners $n_p = 3$ with the corresponding values of $\tau \in \{10^{-8},10^{-4},10^{-2}\}$. This is the same setup as in~\cite{antil2022efficient} and we refer to this paper for additional implementation details. The preconditioners are factorized and their factorizations are stored.

After solving the $j$th set of systems, we have the basis $\bfZ^{(j)}$ with which we can approximate the solutions of $\bfU^{(j)}$; see~\eqref{eqn:snapshot}. Therefore, rather than compute the matrix $\bfU^{(j)}$ and then compress it subsequently, following Remark~\ref{rem:mpgmpres},  we directly compute the new ``snapshot'' matrix as $\widehat\bfS = \bmat{\bfZ^{(1)} & \dots & \bfZ^{(M)}} \in \R^{N_h\times Mn_m}$. This is because the range of $\bfU^{(j)}$ is approximately captured by the range of $\bfZ^{(j)}$, since $\bfZ^{(j)}$ is the search space for the set of shifted systems. Although it is not technically a matrix of snapshots, we still refer to $\widehat\bfS$ as a snapshot matrix. Notice that the number of columns of $\widehat\bfS$ is  much smaller than that of $\bfS$ if $n_m \ll Z$, as is the case in numerical experiments. 

Even with this reduction, the cost of storing and computing the dominant singular vectors can be expensive. To alleviate these costs, we use the single view randomized approach as described in Section~\ref{ssec:rand}. To this end, we consider two sketching matrices $\bfOmega \in \R^{(Mn_m) \times \ell_1}$ 
and $\bfPsi \in \R^{N_h \times \ell_2}$ which are standard Gaussian random matrices. Given the target basis size $K$, we choose $\ell_1 = 2K + 1$ and $\ell_2 = 2\ell_1 + 1.$ We observe that the matrix $\widehat{\bfS}$ can be written as 
\[ \widehat{\bfS} =  \widehat{\bfS}_1 + \dots +  \widehat{\bfS}_M \qquad \widehat{\bfS}_j = \bfZ^{(j)} \bfE_j^\top, \]
where, using MATLAB notation, $\bfE_j = \bfI_{Mn_m}(:, (j-1)*n_m+1:j*n_m)$ for $1 \leq j \leq M$ contains columns from the $(Mn_m) \times (Mn_m)$ identity matrix. Therefore, we can easily update the sketches $\bfY_1$ and $\bfY_2$ for each $j$ and then discard the bases $\bfZ^{(j)}$. In fact, the matrix $\bfOmega$ does not need to be formed all at once; this is beneficial since the basis size $n_m$ may not be known in advance. Consider $\bfOmega = \bmat{\bfOmega_1^\top &  \dots& \bfOmega_{M}^\top}^\top $, where $\bfOmega_j \in \R^{n_m \times \ell_1}$. Then the sketch $\bfY_1$ can be formed as 
\[ \bfY_1 = \bfZ^{(1)}\bfOmega_1 + \dots  + \bfZ^{(M)} \bfOmega_M.\]

The details of this approach are given in Algorithm~\ref{alg:romoffline}.

\paragraph{Online stage} 
In the online stage, the steps involve solving the $K\times K$ systems for $Z = Z_+ + Z_- + 1$, see~\eqref{eqn:projsys} which can be expensive. Rather than factorizing each linear system, we can use a  generalized eigendecomposition. To this end, define $\widehat\bfK = \sum_{j=1}^{n_a}f_t^a(\overline{\bfmu})\widehat{\bfA}$, $\widehat\bfg = \sum_{j=1}^{n_g}f_t^g(\overline{\bfmu}) \widehat\bfg $, and consider a generalized eigendecomposition of the pair $(\widehat\bfK, \widehat\bfM)$; that is, compute a Cholesky factorization $\widehat\bfM = \bfL\bfL^\top$ and compute an eigendecomposition $\bfL^{-1}\widehat\bfK\bfL^\top = \bfU \bfLambda \bfU^\top$. We are guaranteed such an eigendecomposition since both $\bfM$ and $\widehat\bfM$ are symmetric positive definite. Then, we can obtain the solution $\widehat\bfy$ as 
\begin{equation}\label{eqn:projsyseig}\widehat\bfy = (\bfV\bfL^{-\top}\bfU) \left( \sum_{k=-Z_-}^{Z_+} w_k  (\bfLambda + e^{z_k} \bfI)^{-1} \right)\bfU^\top\bfL^{-1}\widehat\bfg. \end{equation}

The modification to Algorithm~\ref{alg:romonline} is straightforward and therefore omitted.

\begin{algorithm}[!ht]
\begin{algorithmic}
\REQUIRE Parameter samples $\mc{M} = \{ \bfmu^{(1)}, \dots, \bfmu^{(M)}\}$, basis size $K$
\STATE  Assemble the matrices $\bfA_t$ for $1 \leq t \leq n_a$ and vectors $\bfg_t$ for $1\leq t \leq n_g$. 
\STATE Compute the nodes $z_k = {k\zeta}$  where $\zeta = 1/\log(1/h)$ and $ -Z_- \leq k \leq Z_+ $ where $Z_+ = \lceil \frac{\pi^2}{2 \zeta^2 }\rceil$, $Z_- = \lceil \frac{\pi^2}{2\zeta^2 }\rceil$
\STATE Draw random matrix $\bfPsi \in \R^{N_h \times \ell_2}$
\STATE Initialize sketches $\bfY_1 = \bfzero \in \R^{N_h\times \ell_1}$, $\bfY_2 = []$
\STATE \COMMENT{Stage 1: Constructing sketches of the snapshot matrix}
\FOR {$j\in \{-Z_-,\dots,Z_+\}$}
\STATE Form matrix $\bfK_j =\sum_{t=1}^{n_a} f_t^a({\bfmu}_j) \bfA_t$ and vector $\bff_j = \sum_{t=1}^{n_g} f_t^g({\bfmu}_j) \bfg_t$
\STATE Solve the systems of linear equations~\eqref{eqn:shiftedsys}
 using MPGMRES-Sh to generate a basis $\bfZ^{(j)}\in \R^{N_h\times n_m}$.
 \STATE Draw random matrix $\bfOmega_j  \in \R^{n_m \times \ell_1}$ 
\STATE Update sketches $\bfY_1 = \bfY_1 + \bfZ^{(j)} \bfOmega_j$ and  $\bfY_2 = [\bfY_2; (\bfZ^{(j)})^\top \bfPsi]$
\ENDFOR
\STATE \COMMENT{Stage 2: Compression of snapshot matrix}
\STATE Compute thin QR $\bfY_1 = \bfQ\bfR$
\STATE Compute $\bfW = (\bfPsi^\top\bfQ)^\dagger\bfY_2^\top$ and its thin SVD $\bfU_W\bfSigma_W\bfV_W^\top$
\STATE Compute $\bfV = \bfQ\bfU_W(:,1:K)$ where $K$ is the basis size
\STATE Precompute $\widehat{\bfA}_t = \bfV^\top \bfA_t \bfV$ for $1 \leq t \leq n_a$, $\widehat{\bfM} = \bfV^\top \bfM \bfV$, and vector $\widehat{\bfg}_t = \bfV^\top \bfM\bfw$
\RETURN Basis vectors $\bfV$, matrices $\{\widehat{\bfA}_t\}_{t=1}^{n_a}$ and vectors $\{\widehat{\bfg}_t \}_{t=1}^{n_g} $
\end{algorithmic}
\caption{Reduced basis approach for fractional elliptic PDE: Offline stage}
\label{alg:romoffline}
\end{algorithm}

\subsection{Computational Cost}\label{ssec:compcosts}
We now summarize the computational cost of the offline and online stages. 

In the offline stage, there are two sources of computational cost: the snapshot matrix construction and compression. The construction essentially involves the cost of applying MPGMRES-Sh to $M$ set of $Z$ shifted systems. Let $m_p$ be the number of preconditioners and $k_m$ be the number of MPGMRES-Sh iterations and their product is $n_m = k_mn_p$; then the cost of one set of MPGMRES-Sh iterations is 
\[ T_{\rm MP} = n_pT_{\rm PC} +  T_{\rm PA} n_m + \mc{O}( N_h Zn_m^2 +  Zn_m^3)) \quad \text{flops},\]
where $T_{\rm PC}$ is the cost of constructing each preconditioner and $T_{\rm PA}$ is the cost of applying the combination of the matrix $\bfC_2$ and the preconditioner $\bfP_j$ (lines 3-4 of Algorithm~\ref{alg:mpgmressh}). The total cost of forming the snapshot matrix is $T_{\rm S} = M T_{\rm MP}$. The cost of the randomized low-rank approximation is analyzed in Section~\ref{ssec:rand}; assuming the basis size $K$, the cost is $T_{\rm R} = \mc{O}(N_h n_m M  K + K^2 (N_h + n_m) + K^3)$ flops. The total cost of the offline stage is then simply the sum of $T_{\rm S}$ and $T_{\rm R}$.

For the offline stage, the na\"ive approach requires us to solve $Z$ systems with $\mc{O}(K^3)$ cost per system, with a total cost of $\mc{O}(ZK^3 + N_hK^2)$ flops. The proposed approach only requires $\mc{O}(K^3)$ cost for the eigendecomposition and $\mc{O}(K^2)$ flops per system solve for a total cost  $\mc{O}(K^3 + K^2 (Z + N_h))$ flops. If $K$ is large, computing the eigendecomposition is computationally expensive or even prohibitive. Then one can use techniques from~\cite{elman2015preconditioning} to further reduce the computational cost. Another option is to view~\eqref{eqn:projsys} through the lens of a generalized Sylvester equation; one can then use techniques from~\cite{simoncini2016computational} to reduce the computational costs.

\section{Numerical Experiments}
In this section, we explore the performance of the proposed ROMs on three illustrative applications. The relative error is measured in the $L^2$ norm with a reference finite element solution computed using MPGMRES-Sh. {The finite element discretization was carried out using continuous piecewise linear elements.} The MPGMRES-Sh solver has been validated in~\cite[Section 5.1]{antil2022efficient}. The timing results were performed on a Mac Mini 2020, with $16$ GB RAM and running MATLAB 2021a. 
\subsection{Gaussian process}\label{ssec:gp}
In this application, we consider simulation from the Gaussian process, by solving the stochastic PDE on the domain $\Omega \in \R^d$
\begin{equation}
    (\kappa^2 - \Delta)^{(\nu + d/2)/2} u(\bfx) = \xi(\bfx),  \qquad \bfx \in \Omega
\end{equation}
where $\xi$ is the white noise process, and in the exponent  $\nu$ is the smoothness parameter and $d$ the spatial dimension. The parameter $\kappa^2$ is related to the correlation length of the process.  We want to compute the solution of the SPDE for multiple values of the parameters $\kappa^2$ and $\alpha = (\nu + d/2)/2$. For this application, there is a single parameter $\bfmu = (\kappa^2)$. Therefore, $n_a = 2$ with $f_1^a(\bfmu) = 1$ and $f_2^a(\bfmu) = \kappa^2$; similarly, $n_g = 1$ with $f_1^g(\bfmu)=  1$. The weak forms are $a_1=\int_\Omega \nabla u\cdot \nabla v d\bfx$, $a_2 = \int_\Omega uv d\bfx$, and $g_1 = \int_\Omega \xi vd\bfx$. 

\begin{figure}[!ht]
    \centering
    \includegraphics[scale=0.4]{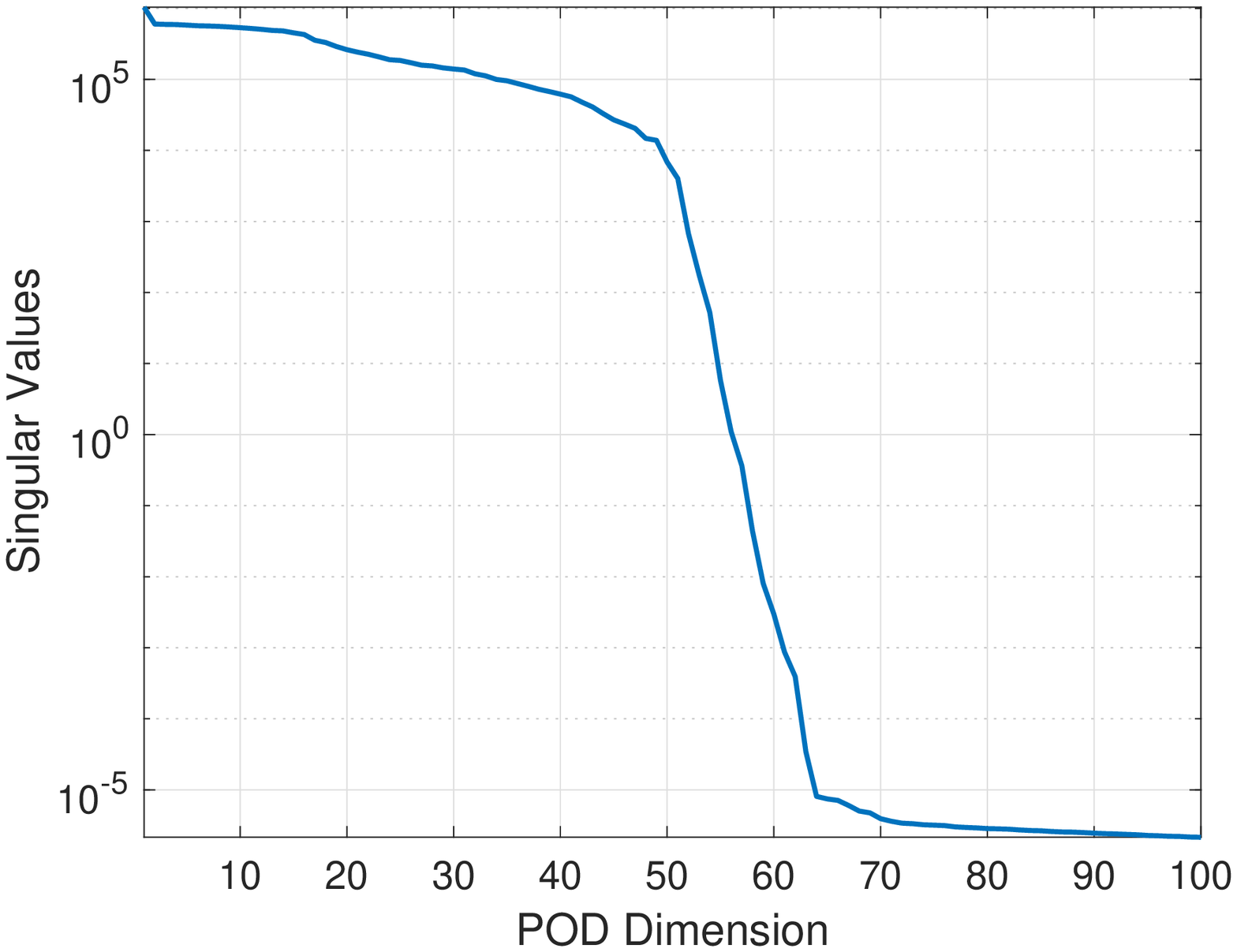}
    \includegraphics[scale=0.4]{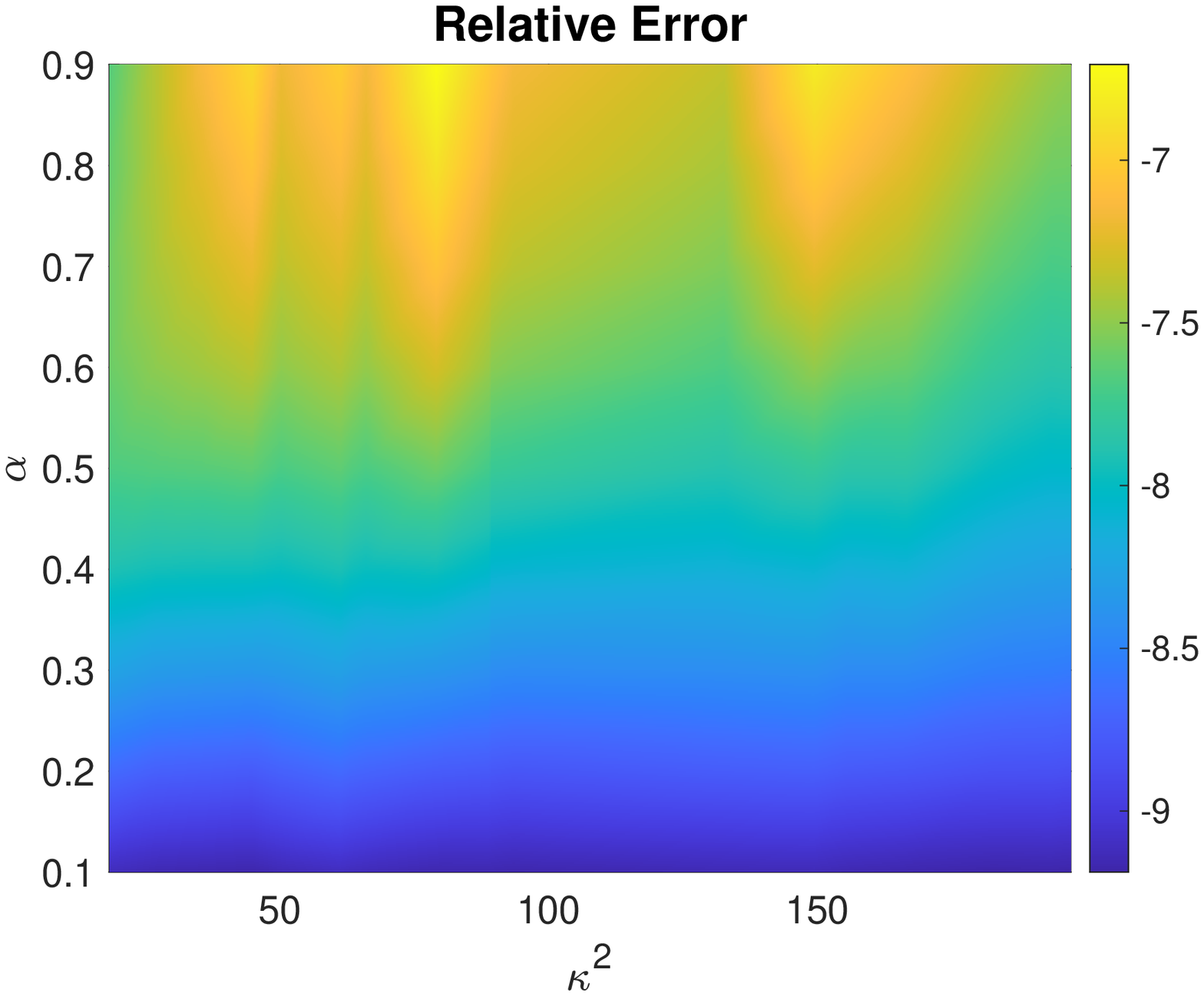}
    \caption{(left) Singular values of the snapshot matrix $\widehat\bfS$. (right) The relative error ($\log_{10}$) in the test set.}
    \label{fig:kappa}
\end{figure}

\paragraph{Experiment 1: Accuracy} We take the domain $\Omega = (0,1)^2$ and consider homogeneous Neumann boundary conditions for the SPDE. The domain is discretized using a grid of size $257\times 257 = 66,049$ points. We use a basis size of $K=100$ based on the decay of the singular values. We generate a training set of parameters $\kappa^2 \in [10,200]$ in increments of $5$, giving a training set of size $|\mc{M}| = 39$. We investigate the accuracy of the approach for several instances of $\kappa^2$ and $\alpha$. Specifically, we generate a test set involving $20$ values of $\kappa^2 \in [10,200]$ and $\alpha \in [0.1,0.9])$ in increments of $0.1$ giving a test set of $180$ candidates. As can be seen from Figure~\ref{fig:kappa}, the relative error  
is small throughout the test set with errors ranging from $10^{-9}-10^{-7}$. The errors show a mild increase with the exponent $\alpha$.

\paragraph{Experiment 2: Timings}In this next experiment, we investigate the reduction in computational costs. We focus on the timing of the compression.  The time of solving~\eqref{eqn:shiftedsys} for a single parameter sample $\bfmu^{(j)}$ using the na\"ive approach is $23.33$ seconds whereas for MPGMRES-Sh the cost is $1.28$ seconds (this time includes the time to factorize the preconditioners).   Furthermore, in the na\"ive approach, the size of the snapshot matrix $\bfS$ is $66,049 \times 11895$, whereas in the proposed approach the size of the snapshot matrix is $66,049 \times 1872$. The size of the MPGMRES-Sh basis for each $\bfmu^{(j)}$  was  $48$. Since the cost of forming and compressing the snapshot matrix is enormous, we do not consider the na\"ive approach for timing purposes anymore. For the baseline approach, we use the approach in Section~\ref{ssec:prop} but use SVD instead of the randomized approach. 

\begin{table}[!ht]
    \centering
    \begin{tabular}{c|c|c}
    & Baseline & Proposed \\ \hline
      Snapshot    &  $47.03$  & $46.66$\\
      Compression   & $18.69$ & $0.59$   
    \end{tabular}
    \caption{Wall-clock time (in seconds) of the time for computing the snapshot matrix $\widehat\bfS$ and the compression stages.}
    \label{tab:compresstime}
\end{table}

We report both the time to construct the snapshot matrices $\widehat\bfS$ and the time for compression (SVD vs randomization). The results are reported in Table~\ref{tab:compresstime}. Both methods give a comparable time for constructing the matrix. However, the baseline approach takes slightly more time (since the total number of MPGMRES-Sh iterations is not known in advance, and memory has to be allocated on-the-fly which makes it inefficient). Note that the proposed approach does not form the matrix $\widehat\bfS$ and the reported time includes the time to form the sketches $\bfY_1$ and $\bfY_2$. However, when it comes to the compression costs, our approach shows a $30\times $ speedup over the SVD. This speedup will be more significant for larger-scale problems. 

Finally, we time for evaluating the ROM for $\alpha = 0.5$ is $0.01$ seconds, which is much more efficient  than both the MPGMRES-Sh {(1.28 seconds)} and the na\"ive approaches {(23.33 seconds)}.
\subsection{Fractional cookies}
\begin{figure}[!ht]
    \centering
    \includegraphics[scale=0.4]{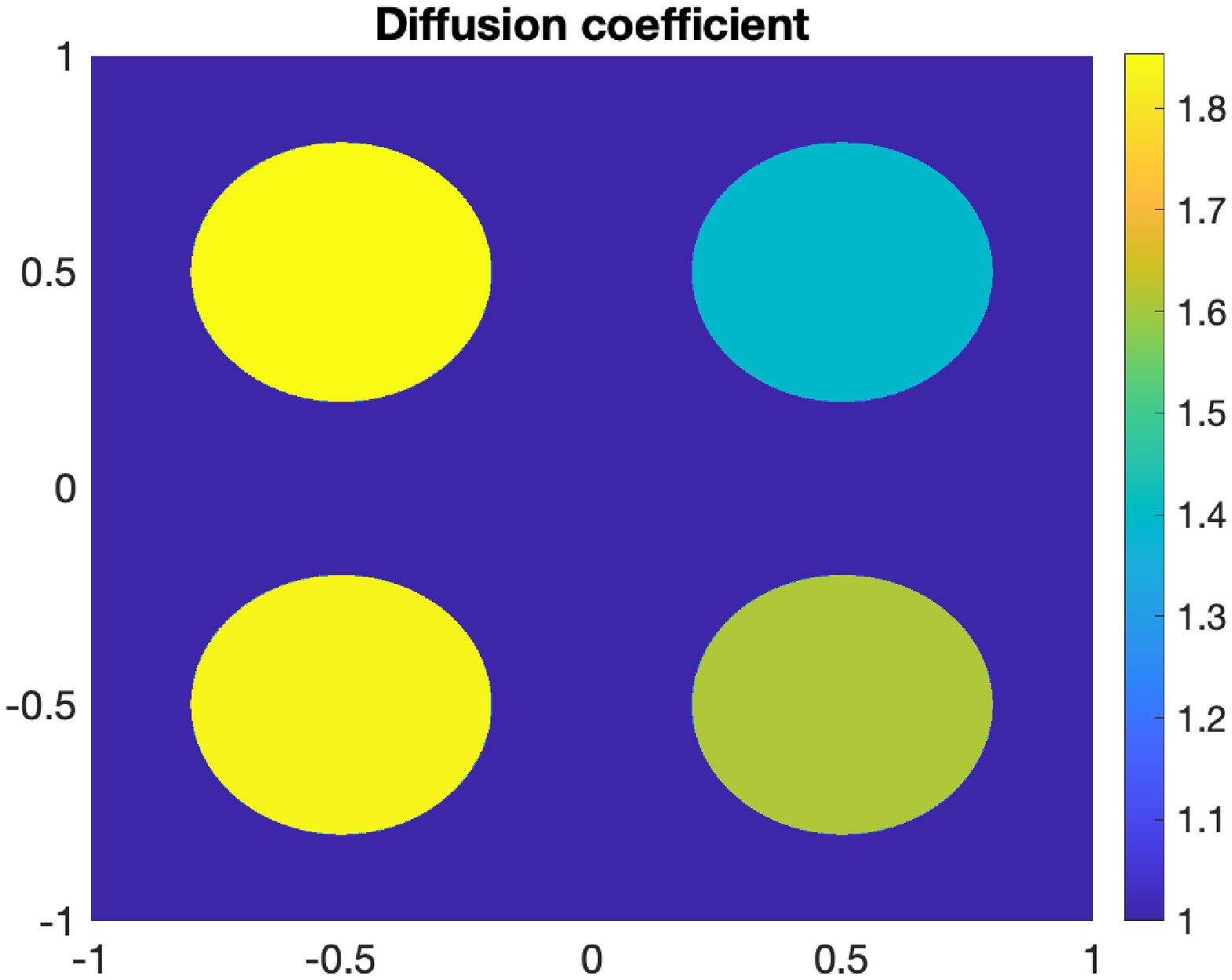}
    \includegraphics[scale=0.4]{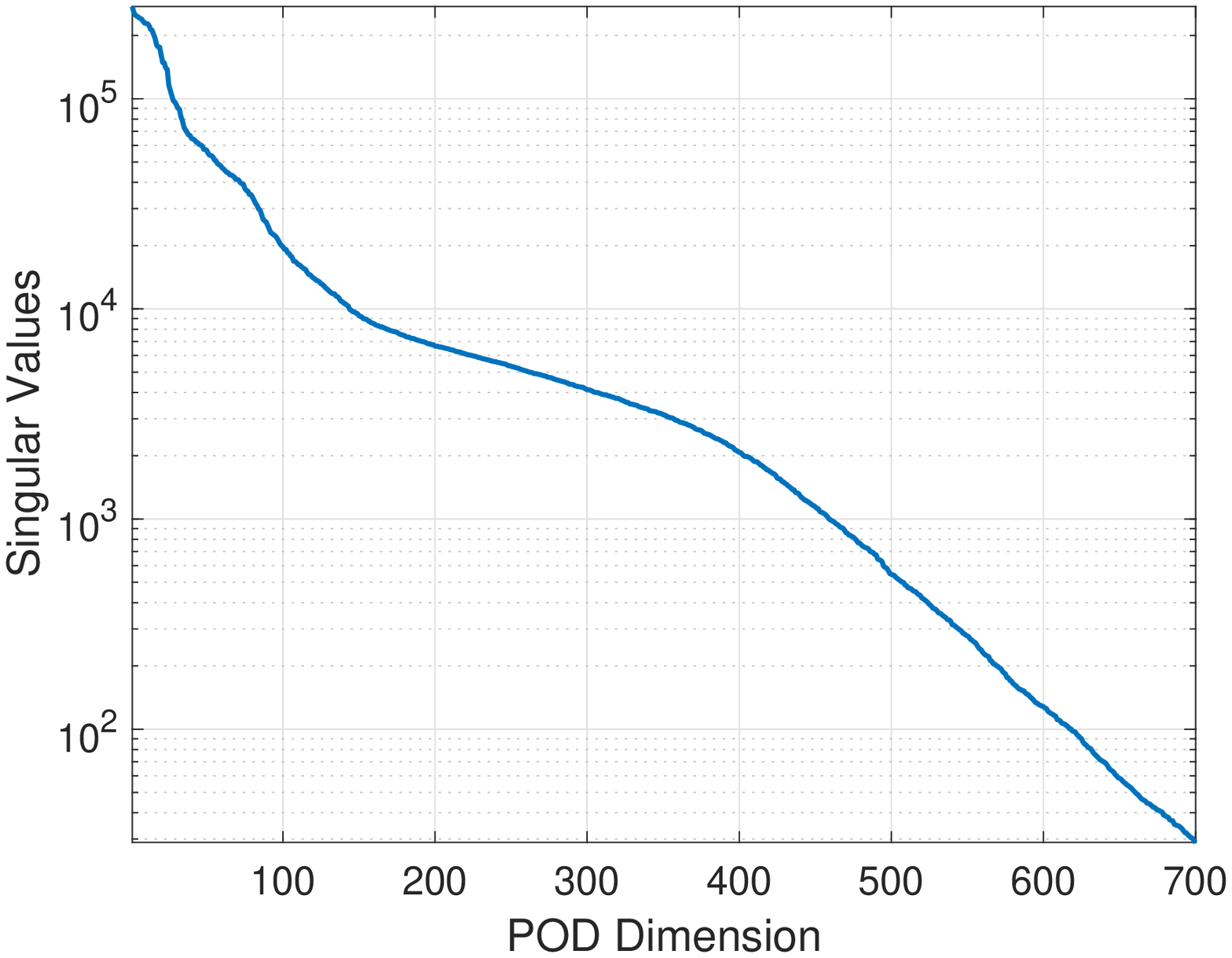}
    \caption{(left) Visualization of the diffusion coefficient $D(\bfx;\bfmu)$ for Case $B$; (right) first $700$ singular values of the snapshot matrix $\widehat\bfS$.}
    \label{fig:fraccookie}
\end{figure}

In this experiment, inspired by the cookies test problem from~\cite{kressner2010krylov}, we consider the fractional parameterized PDE~\eqref{eqn:parampde} with the domain $\Omega = (-1,1)^2$ with Dirichlet boundary conditions,  and the elliptic operator is given by $\mc{A}(\bfx;\bfmu) = -\nabla\cdot(D(\bfx;\bfmu) \nabla (\cdot))$ where 
\[ D(\bfx;\bfmu) = 1 + \sum_{t=1}^p\mu_j\chi_{\Omega_t}(\bfx) \qquad  \]
where $\chi_{A}$ is the characteristic function (that takes the value $1$ if $\bfx \in A \subset \Omega$ and $0$ otherwise), and  $p$ is the number of regions %\todo{I thought $D$ is the coefficient?} 
(``cookies''). We consider two different cases for the discs: 
\begin{figure}[!ht]
    \centering
    \includegraphics[scale=0.4]{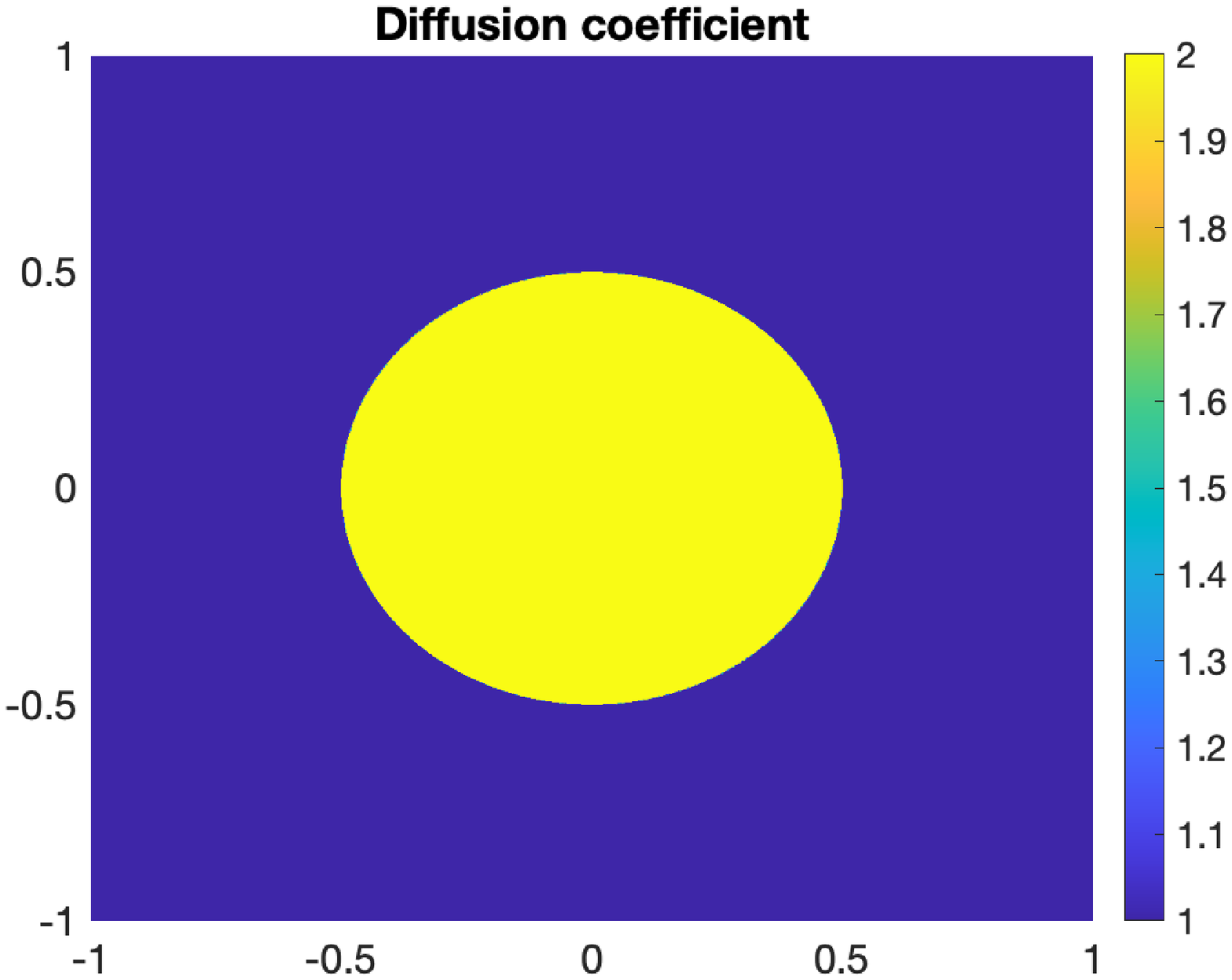}
    \includegraphics[scale=0.4]{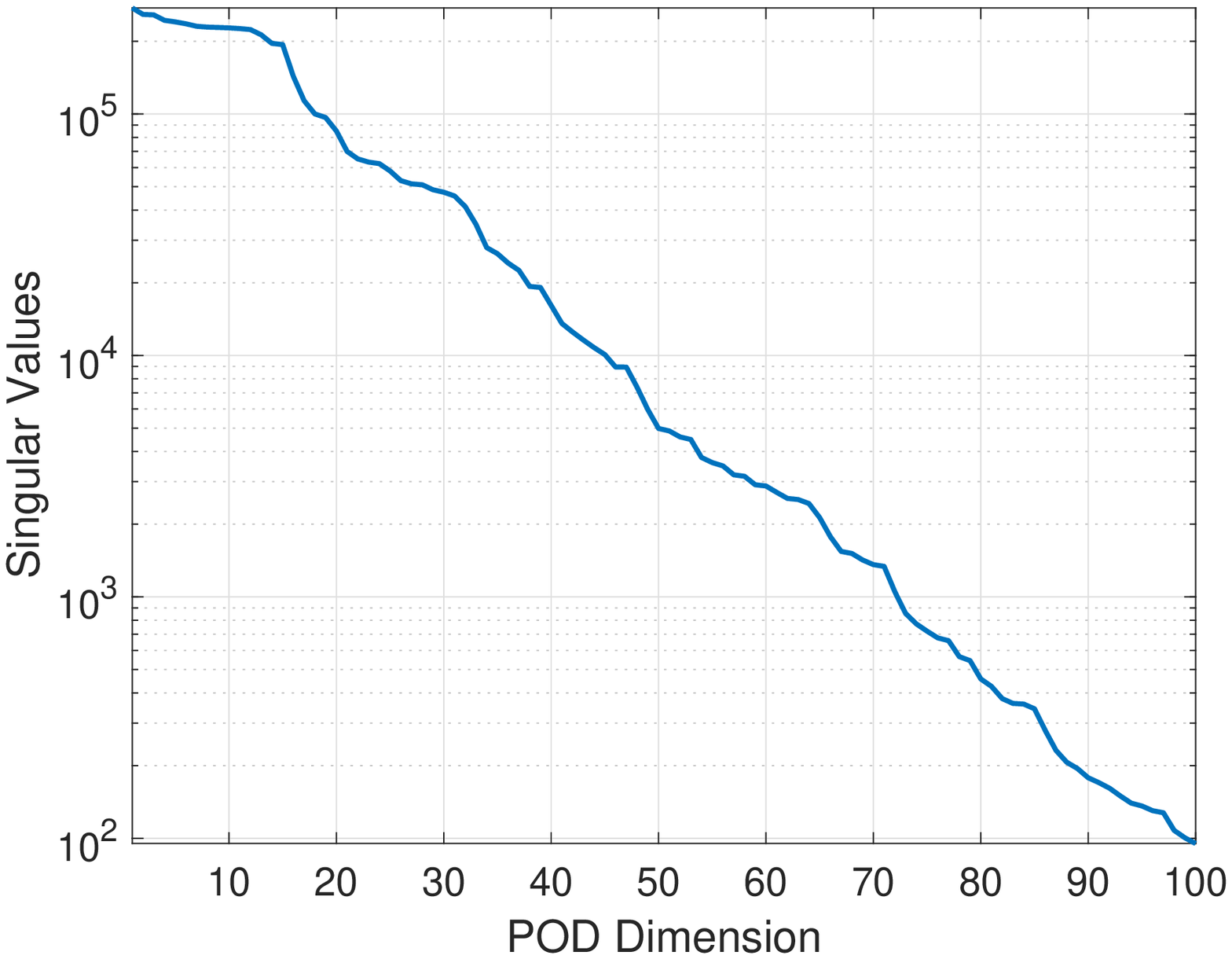}
    \caption{(left) Visualization of the diffusion coefficient $D(\bfx;\bfmu)$; (right) first $100$ singular values of the snapshot matrix $\widehat\bfS$.}
    \label{fig:fraccookie2}
\end{figure}
\begin{enumerate}
    \item Case A: the number of discs $p =1$ and $\Omega_1$ is a disc of radius $0.5$ and center $(0,0)$.
    \item Case B: the number of discs $p=4$ and $\Omega_j$ are discs of radii $0.3$ with centers $(-0.5,-0.5)$, $(-0.5,0.5)$, $(0.5,0.5)$, and $(0.5,-0.5)$.
\end{enumerate} See Figures~\ref{fig:fraccookie2} and~\ref{fig:fraccookie} for visualizations of the diffusion coefficient $D$.

\begin{figure}[!ht]
    \centering
    \includegraphics[scale=0.4]{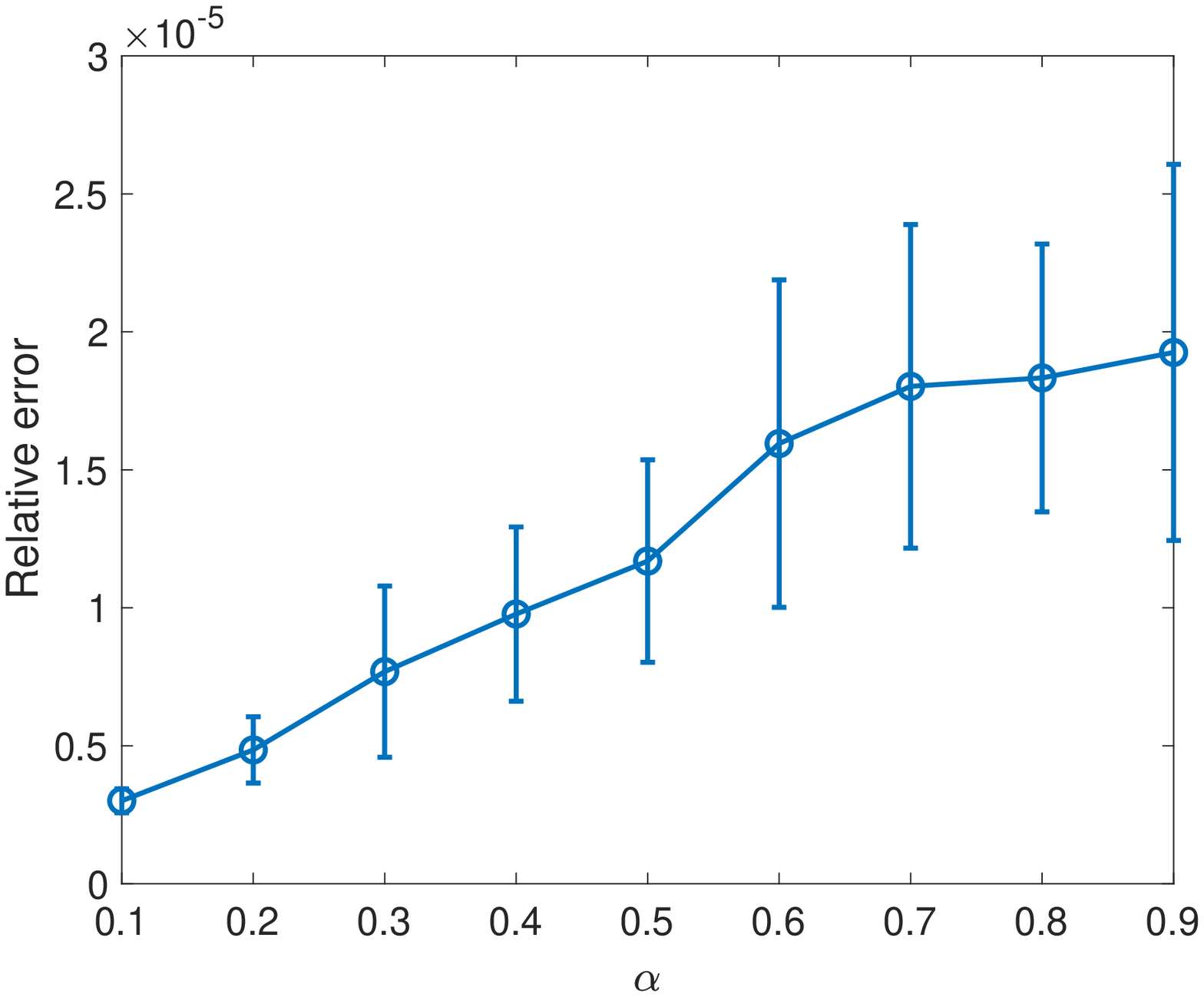}
    \includegraphics[scale=0.4]{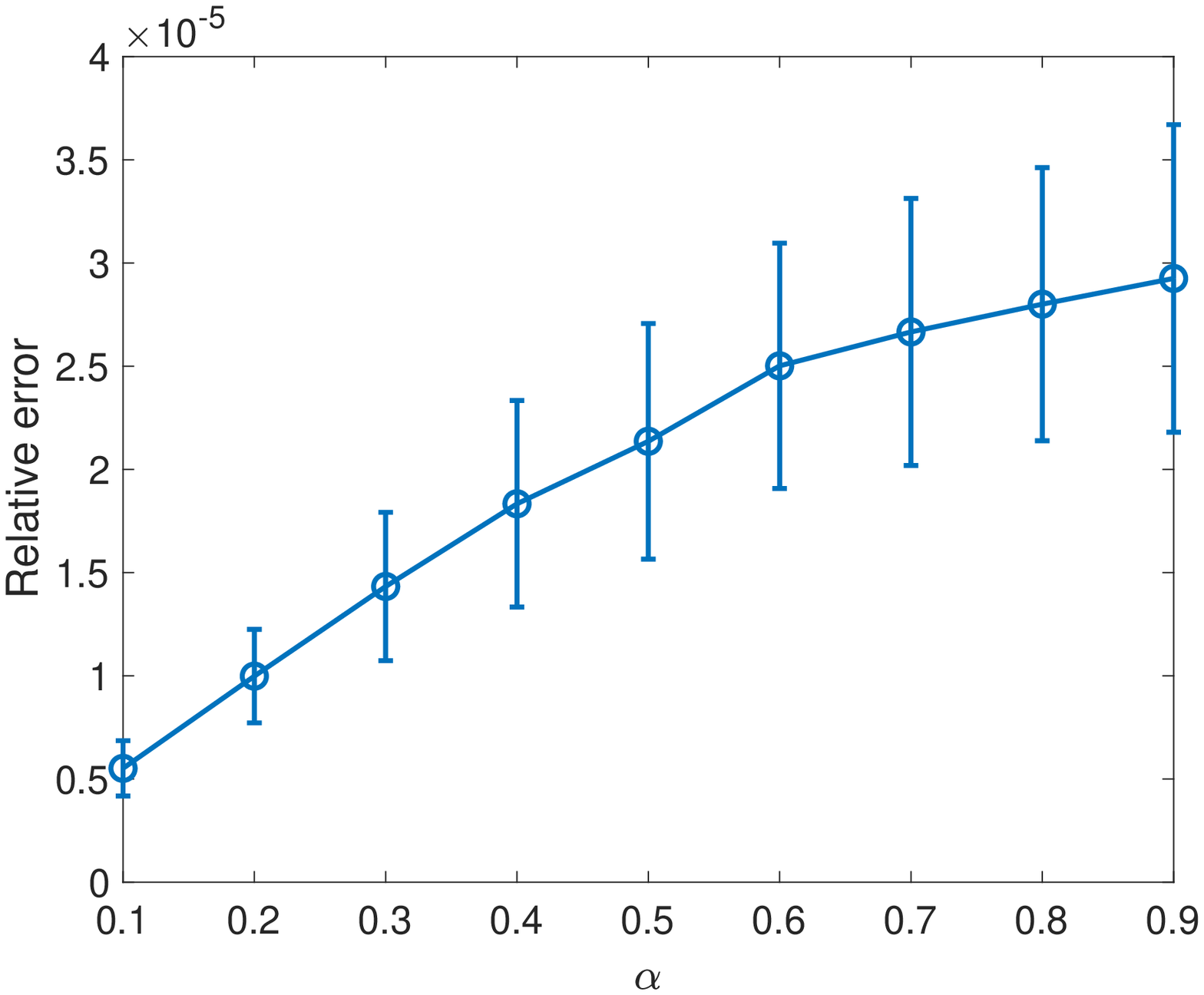}
    \caption{Relative error in the test set (left) Case A, (right) Case B. The mean of the error over $10$ random samples is plotted for different values of $\alpha$, with the error bar denoting one standard deviation.}
    \label{fig:errcookie}
\end{figure}

We take the values of $\bfmu \in [0,1]^p$.  The fractional PDE is augmented with zero Dirichlet boundary conditions. For this problem, $n_a = p+1$ with $a_t(u,v) = \int_{\Omega_t} \nabla u \cdot \nabla v d\bfx$ for $1 \leq t \leq p$ and $a_5 = \int_\Omega \nabla u \cdot \nabla v d\bfx$; similarly, $n_g=1$ and $g_1(v) = \int_\Omega vd\bfx$. 
We discretize the domain using a grid of size $257^2$. To generate the training set, we pick $M=100$ parameters $\bfmu^{(j)}$ from $[0,1]^p$ sampled using Latin hypercube sampling. Based on the decay of the singular values, we take the reduced basis size to be $K =100$ (Case A) and $K = 700$ (Case B). Compared to the test problem in Section~\ref{ssec:gp}, the singular values of the snapshot matrix do not decay as sharply in Case B (see Figure~\ref{fig:fraccookie}).   For the test set, similarly, we consider $10$ parameters uniformly chosen at random from the range $\bfmu \in [0,1]^p$ and for values of $\alpha \in [0.1,0.9]$. {Figure~\ref{fig:errcookie} shows the distribution of relative errors for Case A and Case B.}
In Table~\ref{tab:time_cookies}, we display the computational time for the offline and online stages. As can be seen, for both cases, the online computational time outperforms the Direct and the MPGMRES-Sh algorithms. However, the offline and online times of Case B are higher since the reduced basis has a higher rank $K$. Note that when the reduced basis rank $K$ is high, the ROM may offer computational benefits over the conventional approaches only if $N_h$ is sufficiently high.
\begin{table}[!ht]
    \centering
    \begin{tabular}{c|c|c|c|c|c}
     & Snapshot & Compression & Direct & MPGMRES-Sh & Online \\ \hline
      Case A   & $143.72$  & $0.67$  &  $71.83$ & $1.48$ & $0.006$   \\
     Case B    & $236.49$ & $20.51$ & $71.71$ & $1.93$ & $0.37$
    \end{tabular}
    \caption{Wall-clock time (in seconds) for various components of the cookies test problem. }
    \label{tab:time_cookies}
\end{table}

\subsection{Anisotropic diffusion}\label{ssec:aniso}
In this example, we consider the fractional PDE~\eqref{eqn:paramfpde} with the anisotropic diffusion model $\mc{A}(\bfx;\bfmu) = -\nabla\cdot(\bfTheta(\bfmu) \nabla(\cdot))$, where $\bfTheta$ takes the form 
\[ \bfTheta(\bfmu) = \bmat{ \cos\theta & - \sin\theta\\ \sin \theta & \cos \theta} \bmat{D_1  & \\ & D_2 }  \bmat{ \cos\theta & - \sin\theta\\ \sin \theta & \cos \theta}^\top.\]
The fractional PDE is posed on the domain $\Omega = (0,1)^2$ and is  augmented with zero Dirichlet boundary conditions.  The parameters $\bfmu = (D_1,D_2, \theta)$ include the two diffusion coefficients and the angle $\theta$.

We can write the weak form as $a(u,v; \bfrho) = \sum_{t=1}^3 f_t(\bfmu) a_t(u,v)$, where $a_1(u,v) = \int_{\Omega} \bfe_1 \bfe_1 \nabla u \cdot \nabla v d\bfx$, $a_2 (u,v) = \int_{\Omega} (\bfe_1 \bfe_2^\top + \bfe_2\bfe_1^\top) \nabla u \cdot \nabla v d\bfx$, and $a_3(u,v) = \int_{\Omega} \bfe_2 \bfe_2^\top \nabla u \cdot \nabla v d\bfx$. Here, $\bfe_1,\bfe_2$ are the two columns from the $2\times 2$ identity matrix. Furthermore, $f_1^a(\bfmu) = (\cos^2 \theta) D_1 + (\sin^2 \theta) D_2$,  $f_2^a(\bfmu) = \sin\theta\cos\theta (D_2-D_1)$, and $f_3^a(\bfmu) = (\cos^2 \theta) D_2 + (\sin^2 \theta) D_1$. We take $b(\bfx;\bfmu) = 1$ so that $f_1^g(\bfmu) = 1$ and $g_1(v) = \int_\Omega v d\bfx$. 
\begin{figure}[!ht]
    \centering
    \includegraphics[scale=0.4]{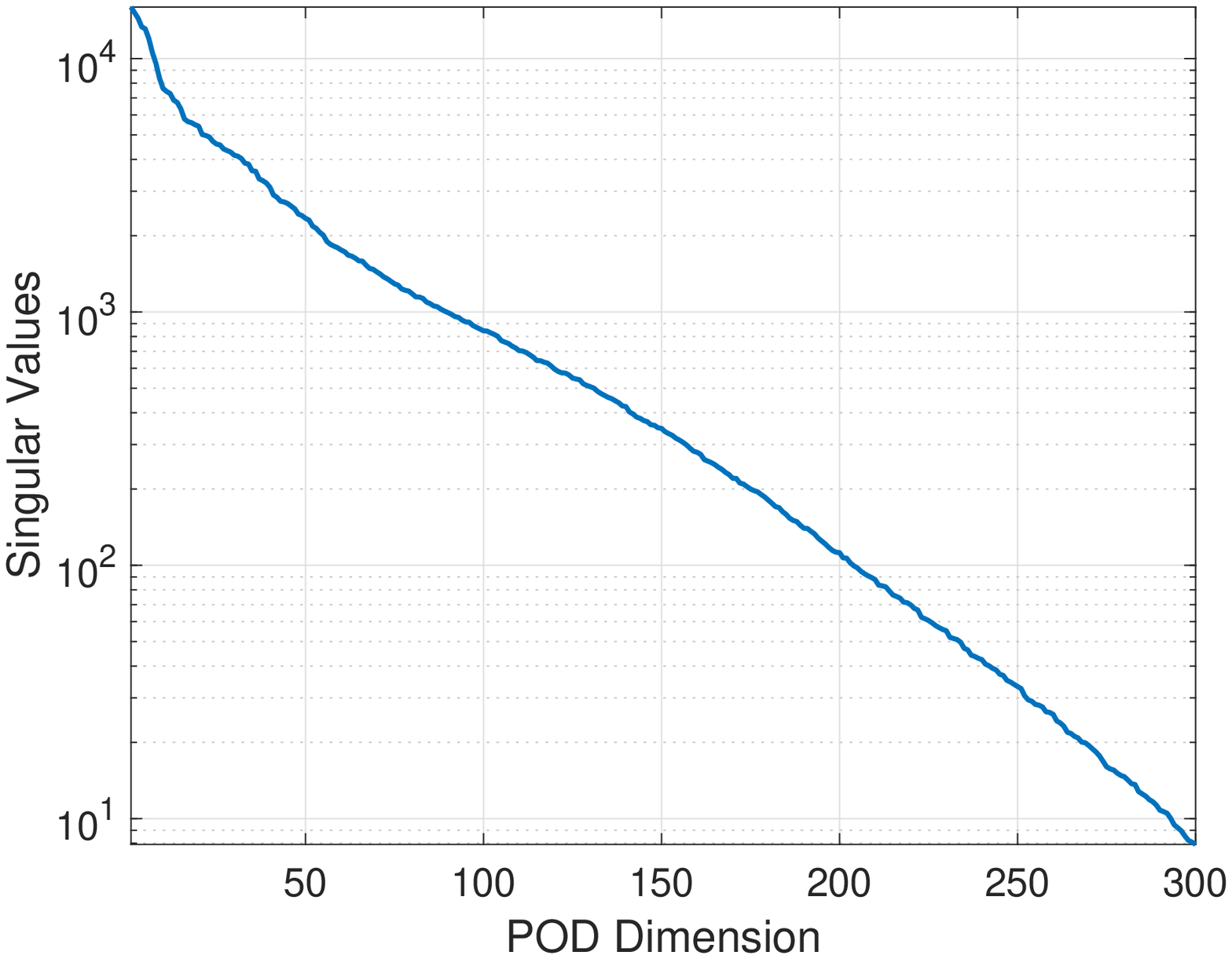}
    \includegraphics[scale=0.4]{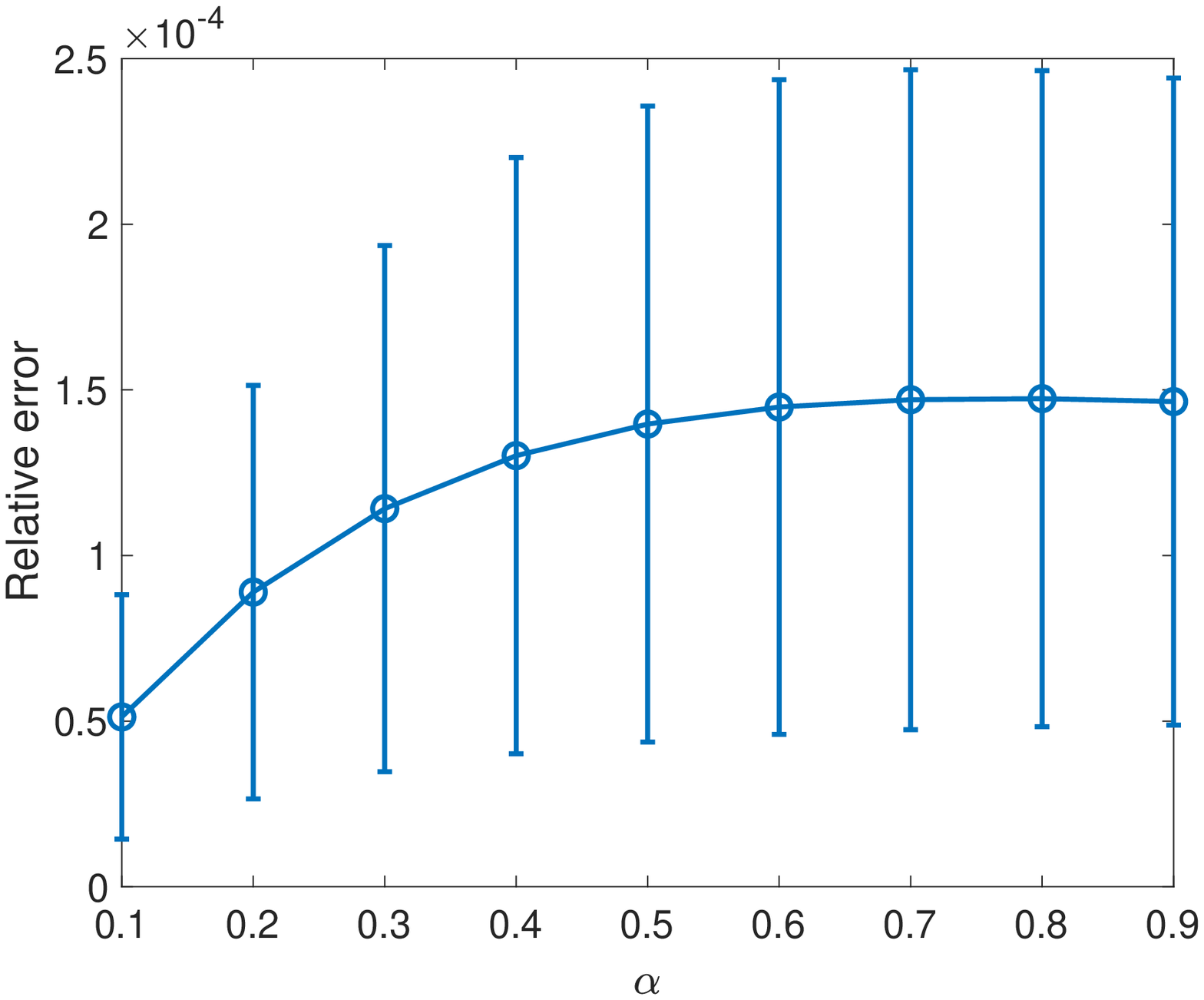}
    \caption{(left) Singular values of the snapshot matrix $\widehat\bfS$, (right) the relative error in the test set. The mean of the error over $10$ random samples is plotted for different values of $\alpha$, with the error bar denoting one standard deviation.  }
    \label{fig:aniso}
\end{figure}

For the training set, we generate $100$ samples generated uniformly at random from the ranges $(D_1, D_2, \theta) \in [0.5,4.5]^2 \times [0,\pi/2]$. The basis size $K$ is taken as $300$. The test set is drawn randomly similarly to the training set and the maximum error over $90$ samples is $\sim 2.5\times 10^{-4}$; the distribution of the error for different values of the exponent $\alpha$ is given in Figure~\ref{fig:aniso}.

\section{Conclusions}
This paper proposes ROMs for parametric fractional elliptic PDEs. The recommended approach combines efficient solvers for shifted linear systems and randomization to dramatically reduce the computational and storage costs over a na\"ive version of the algorithm. Numerical results demonstrate the computational benefits and viability of the ROMs. Future work involves adaptive construction of the ROM to guarantee a certain user-specified error tolerance and the development of error indicators and metrics necessary to develop stopping criteria for the adaptive algorithms. Another line of future inquiry involves adaptively estimating the size of the reduced basis.

\section{Acknowledgements} 
{AKS was supported in part by the National Science Foundation through the grant DMS-1845406.
HA was partially supported by NSF grant DMS-2110263 and the Air Force Office of Scientific 
Research under Award NO: FA9550-22-1-0248. AKS would like to acknowledge Hussam Al Daas for helpful discussions.}
\bibliography{refs}
\bibliographystyle{abbrv}

\end{document}